\documentclass[times]{nmeauth}

\usepackage{moreverb}

\newcommand\BibTeX{{\rmfamily B\kern-.05em \textsc{i\kern-.025em b}\kern-.08em
T\kern-.1667em\lower.7ex\hbox{E}\kern-.125emX}}

\usepackage{url}
\usepackage{graphicx,graphics}
\usepackage{subfigure}
\usepackage{epsfig}
\usepackage{amsmath,amssymb,amsfonts,amsthm}
\usepackage{mathrsfs}
\usepackage{listings}    
\usepackage{color}
\usepackage{float}
\usepackage{tabularx}
\usepackage{multirow}
\usepackage{stmaryrd}

\usepackage{pgfplots}
\usetikzlibrary{plotmarks}

\theoremstyle{remark}

\newcommand{\Eref}[1]{Equation (\ref{#1})}
\newcommand{\fref}[1]{Figure (\ref{#1})}
\newcommand{\Erefs}[1]{Equations (\ref{#1})}
\newcommand{\frefs}[1]{Figures~(\ref{#1})}

\newcommand{\rmd}{\mathrm{d}}

\newcommand{\bigb}{\mathbf{B}}

\newcommand{\dd}{\mathbf{D}}

\newcommand{\kk}{\mathbf{K}}

\newcommand{\uu}{\mathbf{u}}

\newcommand{\xx}{\mathbf{x}}

\def\volumeyear{2013}

\begin{document}

\runningheads{S.~Natarajan, E.T.~Ooi, H.~Man and Ch.~Song}{Strain smoothing and Semi-analytical approach in a quadtree mesh}

\title{{\small SHORT COMMUNICATION} \\
\smallskip
Finite element computations on quadtree meshes: strain smoothing and semi-analytical formulation}

\author{Sundararajan Natarajan\corrauth, Ean Tat Ooi, Hou Man, Chongmin Song}

\address{School of Civil and Environmental Engineering, The University of New South Wales, Sydney, NSW 2052, Australia.}

\corraddr{School of Civil and Environmental Engineering, The University of New South Wales, Sydney, NSW 2052, Australia. Email: sundararajan.natarajan@gmail.com. Tel: +61 9358 5030.}

\begin{abstract}
This short communication discusses two alternate techniques to treat hanging nodes in a quadtree mesh. Both the techniques share similarities, in that, they require only boundary information. Moreover, they do not require an explicit form of the shape functions, unlike the conventional approaches, for example, as in the work of Gupta~\cite{gupta1978} or Tabarraei and Sukumar~\cite{tabarraeisukumar2005}. Hence, no special numerical integration technique is required. One of the techniques relies on the strain projection procedure, whilst the other is based on the scaled boundary finite element method. Numerical examples are presented to demonstrate the accuracy and the convergence properties of the two techniques.
\end{abstract}


\maketitle


\vspace{-6pt}

\section{Introduction}
\vspace{-2pt}

The finite element method (FEM) is widely used approach to solve the partial differential equations (PDEs). The FEM requires the domain to be discretized into non-overlapping regions, called elements. The individual elements are connected together by a topological map called a mesh. When modeling problems involving localized deformation, steep gradients or discontinuous surfaces, a very fine mesh is usually required for accurate results. To reduce the computational time, adaptive refinement techniques are usually preferred over a uniform mesh refinement. Compared with the conforming refinements, quadtree/octree meshes are particularly easy to implement. However, special shape functions have to be constructed on elements with hanging nodes. The shape functions has to be conforming and form a partition of unity.

\subsection{Background} In literature, two approaches were proposed to treat the element with hanging nodes without sub-triangulation. The first approach relied on deriving a set of conforming shape functions for elements with hanging nodes~\cite{gupta1978}. The other approach realised elements with hanging nodes as polygon elements~\cite{tabarraeisukumar2005}. Here, we present a brief overview of the techniques. For more detailed discussion, interested readers are referred to the literature~\cite{gupta1978,tabarraeisukumar2005}.
\subsubsection*{Conforming shape functions} 
\begin{figure}[htpb]
\centering
\scalebox{0.65}{\input{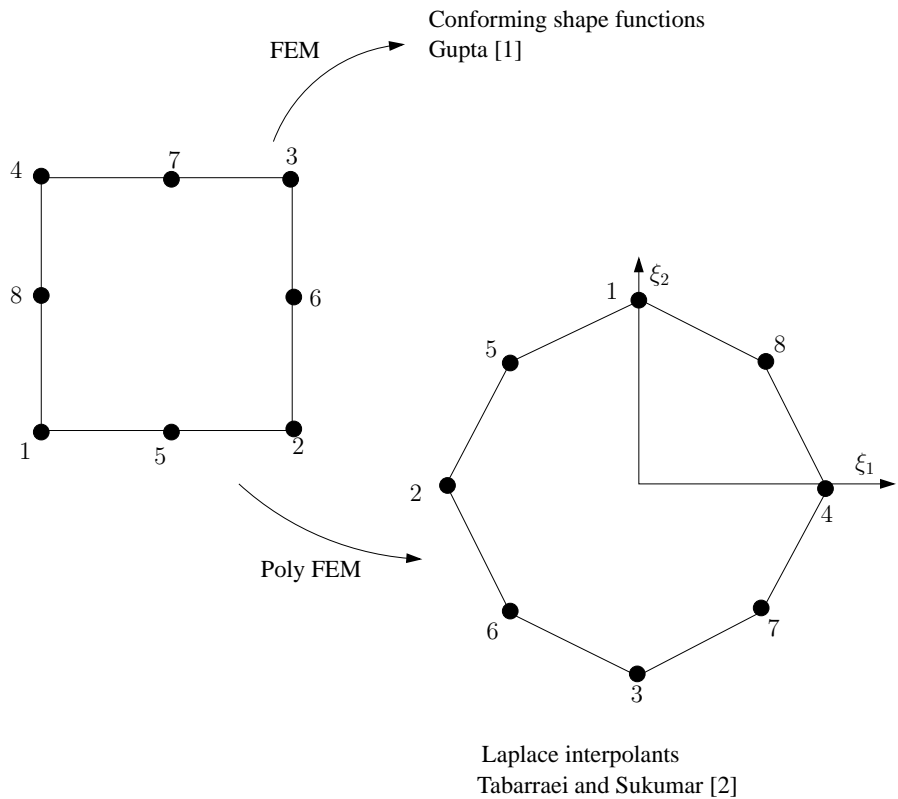}}
\caption{Conventional method to handle hanging nodes in a quadtree mesh: conforming shape functions~\cite{gupta1978} and polygonal FEM~\cite{tabarraeisukumar2005}, where the `filled' circle represents the node.}
\label{fig:hangingnodemethods}
\end{figure}

Gupta~\cite{gupta1978} formulated transition elements with bilinear quadrilateral elements. A conforming set of shape functions for the element with hanging nodes was derived based on the shape functions of the bilinear quadrilateral elements. The shape functions of the corner nodes were modified so that the displacement at the intermediate nodes is equal to the total displacement. The derivatives of the derived shape functions are discontinuous within the element and Gupta~\cite{gupta1978} presented a modified quadrature formula to numerically integrate the terms in the stiffness matrix. The formulation presented was general and can be extended to higher order elements and hexahedral elements as reported in~\cite{lowu2012,mcdillgoldak1987}. The reference shape functions associated with the hanging nodes are given by:
\begin{align}
N^\ast_5(\xi,\eta) &= \frac{1}{2}(1-|\xi|)(1-\eta); \hspace{0.1cm} N^\ast_6(\xi,\eta) = \frac{1}{2}(1+\xi)(1-|\eta|) \nonumber \\
N^\ast_7(\xi,\eta) &= \frac{1}{2}(1-|\xi|)(1+\eta); \hspace{0.1cm} N^\ast_8(\xi,\eta) = \frac{1}{2}(1-\xi)(1-|\eta|)
\end{align}
where $N_i^\ast, (i=5,6,7,8)$ are the shape functions for the mid-side nodes (see \fref{fig:hangingnodemethods}) and are active only if the corresponding hanging node is present along a particular edge. Fries \textit{et al.,}~\cite{friesbyfut2011} employed special elements to associate explicitly the degrees of freedom to the hanging nodes. Though relatively easy to implement in 2D, its extension to 3D is not straightforward, as the hanging nodes can be either on a face or an edge. One possible way to circumvent this difficulty is to restrict the possibilities of various configurations that may arise when building a quad/octree meshes, but this needs to be further investigated. Legrain \textit{et al.,}~\cite{legrainallais2011} treated the hanging nodes by first choosing the right degree of freedoms and then constrain them to ensure continuity of the finite element field.

\subsubsection*{Polygonal FEM} 
\begin{figure}
\centering
\includegraphics[scale=0.5]{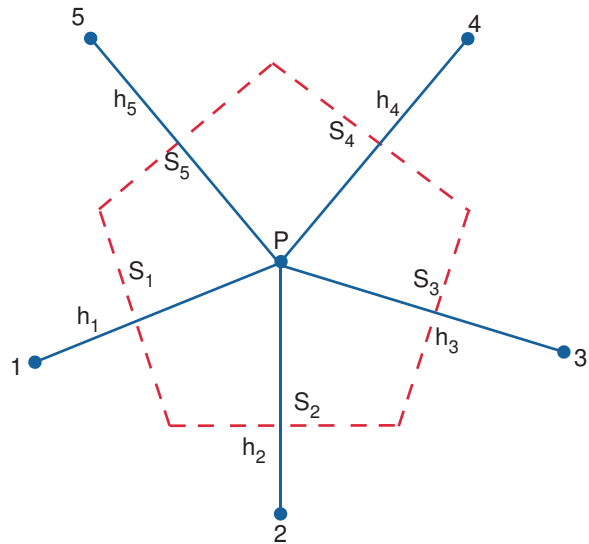}
\caption{Vorono\"{i} diagram of a point $P$.}
\label{fig:vorfig}
\end{figure}

Tabarraei and Sukumar~\cite{tabarraeisukumar2005} considered the elements containing hanging nodes as polygonal elements which requires integration over arbitrary polygons. The numerical integration can be performed by either splitting the elements into several simplices~\cite{sukumarmalsch2006} or by employing a complex mapping~\cite{natarajanbordas2009}. Tabarraei and Sukumar~\cite{tabarraeisukumar2005} employed Laplace interpolants over the polygons. The Laplace interpolant is also called the natural neighbor interpolant~\cite{tabarraeisukumar2005}. It provides a natural weighting function for irregularly spaced nodes. For a point $P$ with $n$ natural neighbors, the Laplace shape functions for node $P_I$ can be written as:
\begin{equation}
\phi_I(\xx) = \frac{ \alpha_I(\xx)}{\sum\limits_{I=1}^n \alpha_J(\xx)}, \hspace{0.25cm} \alpha_J(\xx) = \frac{s_J(\xx)}{h_J(\xx)}, \hspace{0.2cm} \xx \in \mathbb{R}^2
\end{equation}
where $\alpha_I(\xx)$ is the Laplace weight function, $s_I(\xx)$ is the length of the Vorono\"{i} edge associated with $P$ and node $P_I$, and $h_I(\xx)$ is the Euclidean distance between $P$ and $P_I$ (see \fref{fig:vorfig}).

\subsection{Approach}
In this paper, we present two alternate approaches to treat hanging nodes in a quadtree mesh. The spirit of the approaches proposed here shares some similarities with the work of Tabarraei and Sukumar~\cite{tabarraeisukumar2005}, but instead of using polygonal elements, the stiffness matrix of the elements with hanging nodes are computed either by the strain smoothing technique or by the scaled boundary formulation. The proposed method does not require special numerical integration techniques to compute the stiffness matrix. Both the techniques require only boundary information. An explicit form of shape functions is not required. 

Compared with the standard error estimators~\cite{zienkiewiczzhu1992,tabbarablacker1994,gonzalez-estradanatarajan2012}, the refinement criteria employed in this study is inquisitive. The mesh is refined where steep gradients are expected or near interfaces. The main objective of this short communication is to disseminate the idea of using either the cell-based smoothing technique or the scaled boundary polygon formulation for quadtree/octree meshes. We restrict ourselves to quadtree meshes, although extension to octree meshes and coupling with accurate error estimators~\cite{zienkiewiczzhu1992,tabbarablacker1994,gonzalez-estradanatarajan2012} for adaptive refinement is feasible.

\subsection{Outline}
The paper is organised as follows. Section \ref{sfemsbfem} presents the idea of the cell-based strain smoothing method and the scaled boundary formulation as applied to quadtree meshes. The effectiveness and the robustness of the approaches are demonstrated with a few of benchmark problems in the context of scalar fields and an application to the partition of unity method. Concluding remarks are presented in the last section.

\vspace{-6pt}

\section{Treating hanging nodes based on smoothing and semi-analytical method}
\label{sfemsbfem}
\subsection{Cell-based strain smoothing} 
\begin{figure}[htpb]
\centering
\scalebox{0.65}{\input{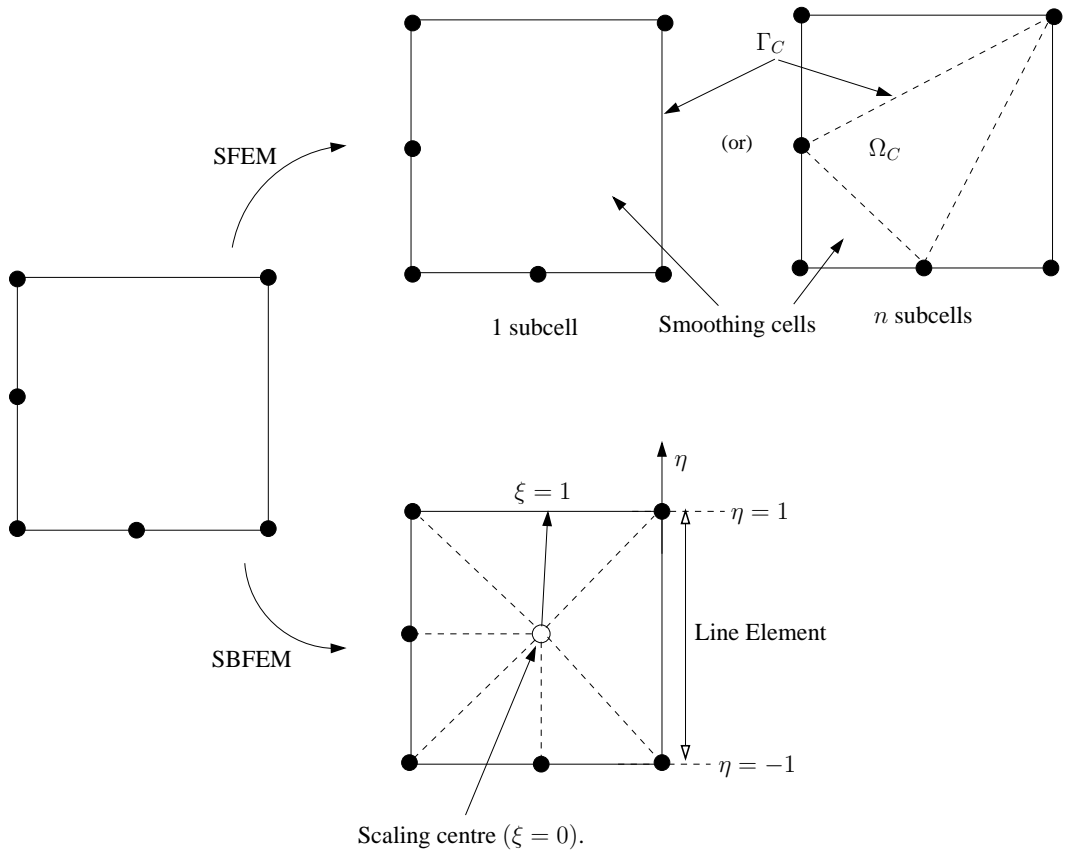}}
\caption{Proposed alternate method to handle hanging nodes in a quadtree mesh: Strain smoothing and the scaled boundary polygon formulation, where the `filled' circle represents the node.}
\label{fig:hangingnodemethods1}
\end{figure}

In this section, the gradient smoothing method, esp., the cell based smoothing method (CSFEM), which is based on the work of Chen and Wang~\cite{chenwang2000}, is briefly discussed. In CSFEM, the elements are divided into smoothing cells over which the standard strain field is smoothed (see \fref{fig:hangingnodemethods1}). This smoothing enables the volume integration of the stiffness matrix to be transformed into a surface integral using the divergence theorem. One way to create subcell geometries is to decompose directly the elements of an existing mesh, the simplest decomposition being realized where each subcell coincides with an element. Other node, edge and face-based decomposition techniques were proposed in a number of papers. For more detailed discussion and derivation for various smoothed FEMs, interested readers are referred to the literature~\cite{liutrung2010,bordasnatarajan2011} and the references therein. 
The strain field, $\tilde{\varepsilon}_{ij}^h$, used to compute the stiffness matrix is computed by a weighted average 	of the standard strain field $\varepsilon_{ij}^h$. At a point $\xx_{C}$ in an element $\Omega^h$, the smoothed strain field is given by:
\begin{equation}
\tilde{\varepsilon}_{ij}^h = \int\limits_{\Omega^h} \varepsilon_{ij}^h (\xx) \Phi(\xx - \xx_C)~\mathrm{d}\xx
\end{equation}
where $\Phi$ is a smoothing function that generally satisfies the following properties:
\begin{equation}
\Phi \geq 0 \hspace{0.5cm} \textup{and} \hspace{0.5cm} \int\limits_{\Omega^h} \Phi(\xx)~\mathrm{d}\xx = 1
\end{equation}
One possible choice of $\Phi$ is given by
\begin{equation*}
\renewcommand{\arraystretch}{1.5}
\Phi = \left\{ \begin{array}{cc} \frac{1}{A_C} & \xx_C \in \Omega_C \\ 0 & \xx_C \notin \Omega_C \end{array} \right.
\end{equation*}
where $A_C$ is the area of the subcell. The CSFEM can be recast within a Hellinger-Reissner variational principle, where the assumed strain is the (constant) smoothed strain over each subcell. For stability, the approximation of the smoothed strain and of the displacement field must satisfy the LBB condition~\cite{babuvska1973}. The smoothed element stiffness matrix for an element $e$ is computed by the sum of the contributions of the subcells
\begin{equation}
\tilde{\kk}^e = \sum\limits_{C=1}^{nC} \int_{\Omega_C} \tilde{\bigb}_C^{\rm T} \dd \tilde{\bigb}_C ~\rmd \Omega
\end{equation}
where $nC$ is the number of smoothing cells of the element. In this study, for the elements with hanging nodes, one subcell and $n$ triangular subcells are considered. The element with hanging nodes is triangulated solely for the purpose of the numerical integration and no additional degree of freedom is introduced. For standard elements, 4 quadrilateral subcells are used, based on the work of Nguyen-Xuan \textit{et al.,}~\cite{nguyen-xuanbordas2008}.

\subsection{Semi-analytical approach: scaled boundary polygon formulation}


\fref{fig:hangingnodemethods1} shows an element with hanging nodes modelled by the SBFEM. The SBFEM~\cite{wolfsong2000} is a semi-analytical computational technique that reduces the governing partial differential equations to a set of ordinary differential equations. In the SBFEM, a local coordinate system is introduced ($(\xi,\eta)$, see \fref{fig:hangingnodemethods1}), where, with reference to \fref{fig:hangingnodemethods1}, $\xi$ is the radial coordinate that with $\xi=0$ at the scaling centre and $\xi=1$ at the cell boundary and $\eta$ is the local coordinate of the one dimensional finite elements discretizing the boundaries of the cells. A scaling centre $O$ is selected at a point from which the whole boundary of the domain is visible. The displacement field covered by a line sector on the boundary of a quadtree cell is approximated by
\begin{align}
\mathbf{u}(\xi,\eta)= & \mathbf{N}(\eta)\mathbf{u}(\xi)\label{eq:dispapprox}
\end{align}
Using the standard strain-displacement relations, the scaled boundary transformation of geometry and \Eref{eq:dispapprox}, the strain field $\boldsymbol{\varepsilon}(\xi,\eta)$ is approximated by
\begin{align}
\boldsymbol{\varepsilon}(\xi,\eta)= & -\mathbf{B}_{1}(\eta)\mathbf{u}(\xi),_{\xi}+\xi^{-1}\mathbf{B}_{2}(\eta)\mathbf{u}(\xi)\label{eq:strainraw}
\end{align}
where $\mathbf{B}_{1}(\eta)$ and $\mathbf{B}_{2}(\eta)$ are the scaled boundary strain displacement matrices. The equilibrium condition for a polygon can be formulated from the principle of virtual work
\begin{align}
\int_{\Omega}\delta\boldsymbol{\varepsilon}^{\mathrm{T}}\boldsymbol{\sigma}d\Omega & =\int_{\Gamma}\delta\mathbf{u}^{\mathrm{T}}\mathbf{t}d\Gamma+\int_{\Omega}\delta\mathbf{u}^{\mathrm{T}}\mathbf{b}d\Omega\label{eq:vitwork}
\end{align}
where $\boldsymbol{\sigma}(\xi,\eta)$ is the stress field, $\delta\boldsymbol{\varepsilon}(\xi,\eta)$ is the virtual strain field, $\delta\mathbf{u}(\xi,\eta)$ is the virtual displacement field, $\mathbf{b}$ is the body force intensity and $\mathbf{t}$ is the surface traction. 

The case with zero body forces $\mathbf{b}=0$ is considered first. Substituting \Erefs{eq:dispapprox} and (\ref{eq:strainraw}), and using the Hooke's law $\boldsymbol{\sigma}=\mathbf{D}\boldsymbol{\varepsilon}$ with $\mathbf{D}$ being the elastic constitutive matrix into \Eref{eq:vitwork} results in~\cite{deekswolf2002}:
\begin{align}
\left.\delta\mathbf{u}_{\mathrm{b}}^{\mathrm{T}}\left(\mathbf{E}_{0}\mathbf{u}(\xi),_{\xi}+\mathbf{E}_{1}^{\mathrm{T}}\mathbf{u}(\xi)-\mathbf{p}\right)\right|_{\xi=1} \nonumber \\
-\delta\mathbf{u}(\xi)^{\mathrm{T}}\int_{0}^{1}\left(\mathbf{E}_{0}\xi\mathbf{u}(\xi),_{\xi\xi}+(\mathbf{E}_{0}+\mathbf{E}_{1}^{\mathrm{T}}-\mathbf{E}_{1})\mathbf{u}(\xi),_{\xi}-\xi^{-1}\mathbf{E}_{2}\mathbf{u}(\xi)\right)d\xi= & 0\label{eq:virtualworkexp}
\end{align}
where $\mathbf{p}$ is the equivalent nodal forces due to the distributed load on the boundary $\mathbf{f}_{\mathrm{t}}$ and
$\mathbf{u}_{\mathrm{b}}$ is the nodal displacement vector on the boundary $(\xi=1)$. $\mathbf{E}_{0}$, $\mathbf{E}_{1}$ and $\mathbf{E}_{2}$ are coefficient matrices that depend only on the geometry and boundary conditions~\cite{deekswolf2002}. \Eref{eq:virtualworkexp} is satisfied if the following conditions
are met
\begin{align}
\delta\mathbf{u}_{\mathrm{b}}^{\mathrm{T}}\left(\mathbf{E}_{0}\mathbf{u}(\xi),_{\xi}+\mathbf{E}_{1}^{\mathrm{T}}\mathbf{u}(\xi)-\mathbf{p}\right)= & 0\label{eq:cond1}\\
\delta\mathbf{u}(\xi)^{\mathrm{T}}\int_{0}^{1}\left(\mathbf{E}_{0}\xi\mathbf{u}(\xi),_{\xi\xi}+(\mathbf{E}_{0}+\mathbf{E}_{1}^{\mathrm{T}}-\mathbf{E}_{1})\mathbf{u}(\xi),_{\xi}-\xi^{-1}\mathbf{E}_{2}\mathbf{u}(\xi)\right)d\xi= & 0\label{eq:cond2}
\end{align}
Invoking the arbitrariness of the virtual displacement functions $\delta\mathbf{u}(\xi)$ in \Eref{eq:cond2} leads to the scaled boundary finite element equation in displacement~\cite{deekswolf2002}
\begin{align}
\mathbf{E}_{0}\xi^{2}\mathbf{u}(\xi)_{,\xi\xi}+(\mathbf{E}_{0}-\mathbf{E}_{1}+\mathbf{E}_{1}^{\mathrm{T}})\xi\mathbf{u}(\xi)_{,\xi}-\mathbf{E}_{2}\mathbf{u}(\xi)= & 0\label{eq:sbfedisp}
\end{align}
which can be solved to obtain $\mathbf{u}(\xi)$ as
\begin{align}
\mathbf{u}(\xi)= & \boldsymbol{\Phi}_{\mathrm{u}}\xi^{-\boldsymbol{\Lambda}_{\mathrm{n}}}\mathbf{c}_{\mathrm{n}}\label{eq:uksi}
\end{align}
where $\boldsymbol{\Lambda}_{\mathrm{n}}$ represent the
eigenvalue matrix with real parts satisfying $\mathrm{Re}(\lambda(\boldsymbol{\Lambda}_{\mathrm{n}}))<0$
and $\boldsymbol{\Phi}_{\mathrm{u}}$ are the eigenvectors corresponding
to the modal displacements in an element hat are obtained from an
eigenvalue decomposition of the Hamiltonian matrix $\mathbf{Z}$
\begin{align}
\mathbf{Z}= & \left[\begin{array}{cc}
\mathbf{E}_{0}^{-1}\mathbf{E}_{1}^{\mathrm{T}} & -\mathbf{E}_{0}^{-1}\\
\mathbf{E}_{1}\mathbf{E}_{0}^{-1}\mathbf{E}_{1}^{\mathrm{T}}-\mathbf{E}_{2} & -\mathbf{E}_{1}\mathbf{E}_{0}^{-1}
\end{array}\right]\label{eq:Hamiltonian}
\end{align}
The integration constants $\mathbf{c}_{\mathrm{n}}$ are determined from the nodal displacements $\mathbf{u}_{\mathrm{b}}=\mathbf{u}(\xi=1)$ as
\begin{align}
\mathbf{c}_{\mathrm{n}}= & \boldsymbol{\Phi}_{\mathrm{u}}^{-1}\mathbf{u}_{\mathrm{b}}\label{eq:intconst}
\end{align}
Substituting \Eref{eq:uksi} into \Eref{eq:cond1} and using \Eref{eq:intconst} results in
\begin{align}
\mathbf{K}_{\mathrm{p}}\mathbf{u}_{\mathrm{b}}= & \mathbf{p}\label{eq:kuf}
\end{align}
where
\begin{align}
\mathbf{K}_{\mathrm{p}}= & \mathbf{E}_{0}\boldsymbol{\Phi}_{\mathrm{u}}\boldsymbol{\Lambda}_{\mathrm{n}}\boldsymbol{\Phi}_{\mathrm{u}}^{-1}+\mathbf{E}_{1}^{\mathrm{T}}\label{eq:stiff}
\end{align}
In the case of non-zero body forces $\mathbf{b}\neq0$, the equivalent nodal forces due to the body load can be expressed in terms of the nodal displacements $\mathbf{u}_{\mathrm{b}}$ by substituting first,~\Eref{eq:uksi} into \Eref{eq:dispapprox} and then into the second term on the right-hand-side of \Eref{eq:vitwork}, resulting in
\begin{align}
\int_{\Omega}\delta\mathbf{u}^{\mathrm{T}}\mathbf{b}d\Omega= & \delta\mathbf{u}_{\mathrm{b}}^{\mathrm{T}}\int_{0}^{1}\boldsymbol{\Phi}_{\mathrm{u}}^{-\mathrm{T}}\xi^{-\boldsymbol{\Lambda}_{\mathrm{n}}+\mathbf{I}}\int_{-1}^{1}\boldsymbol{\Phi}_{\mathrm{u}}^{\mathrm{T}}\mathbf{N}(\eta)^{\mathrm{T}}\mathbf{b}|\mathbf{J}(\eta)|d\eta d\xi\label{eq:body1}
\end{align}
where $|\mathbf{J}(\eta)|$ is the determinant of the Jacobian on the boundary required for coordinate transformation. For many problems, the body force intensity $\mathbf{b}$ can be expressed locally in each quadtree cell as a power function in $\xi$ as
\begin{align}
\mathbf{b}(\xi,\eta)= & \xi^{k}\mathbf{b}(\eta)\label{eq:b_ksi_eta}
\end{align}
Substituting \Eref{eq:b_ksi_eta} into \Eref{eq:body1} and integrating analytically in $\xi$ results in
\begin{align}
\delta\mathbf{u}_{\mathrm{b}}^{\mathrm{T}}\mathbf{p}_{\mathrm{b}}= & \delta\mathbf{u}_{\mathrm{b}}^{\mathrm{T}}\boldsymbol{\Phi}_{\mathrm{u}}^{-\mathrm{T}}\left(-\boldsymbol{\Lambda}_{\mathrm{n}}+(k+2)\mathbf{I}\right)^{-1}\boldsymbol{\Phi}_{\mathrm{u}}^{\mathrm{T}}\int_{-1}^{1}\mathbf{N}(\eta)^{\mathrm{T}}\mathbf{b}(\eta)|\mathbf{J}(\eta)|d\eta d\xi\label{eq:body2}
\end{align}
where $\mathbf{p}_{\mathrm{b}}$ is the equivalent load vector due to the body forces $\mathbf{b}(\xi,\eta)$. \Eref{eq:kuf} can therefore be rewritten as
\begin{align}
\mathbf{K}_{\mathrm{p}}\mathbf{u}_{\mathrm{b}}= & \mathbf{p}+\mathbf{p}_{\mathrm{b}}
\end{align}




\section{Numerical Examples} \label{numexample}
In this section, we study the convergence and the accuracy of the strain smoothing and the SBFEM over quadtree meshes by solving the Laplace and the Poisson's equation over a square domain with appropriate Dirichlet boundary conditions. The numerical results are compared with the analytical solution and with the conforming shape functions and Laplace interpolants. In the last part of the section, the present approaches are combined with the partition of unity method, esp, the extended finite element method (XFEM) to solve a problem involving weak discontinuity. We employ the following convention when presenting the numerical results:
\begin{itemize}
\item FEM - conventional FEM with conforming shape functions~\cite{gupta1978}. A modified Gaussian quadrature is adopted as given in~\cite{gupta1978}.
\item PFEM - polygonal FEM with Laplace interpolants. For the purpose of numerical integration, sub-triangulation with sixth order Dunavant quadrature rule over each sub-triangle is employed.
\item nSFEM - $n$-sided smoothed finite element method. Along each side of the element, the shape functions are assumed to be linear and hence, only one Gauss point is used to numerically integrate the terms in the stiffness matrix.
\item SBFEM - scaled boundary polygon formulation. Only the boundary of the element is discretised with 1D linear shape functions and one point Gauss quadrature is adopted.
\end{itemize}
For the purpose of error estimation, the relative error in $L^2$ norm is used and is given by:
\begin{equation}
|| \uu - \uu^h||_{L^2(\Omega)} = \sqrt{ \int_\Omega [(\uu-\uu^h) \cdot (\uu-\uu^h)]~\rmd \Omega}
\end{equation}
where $\uu^h$ is the numerical solution and $\uu$ is the analytical solution or a reference solution.

\subsection{Patch test} 
\begin{figure}[htpb]
\centering
\includegraphics[scale=0.8]{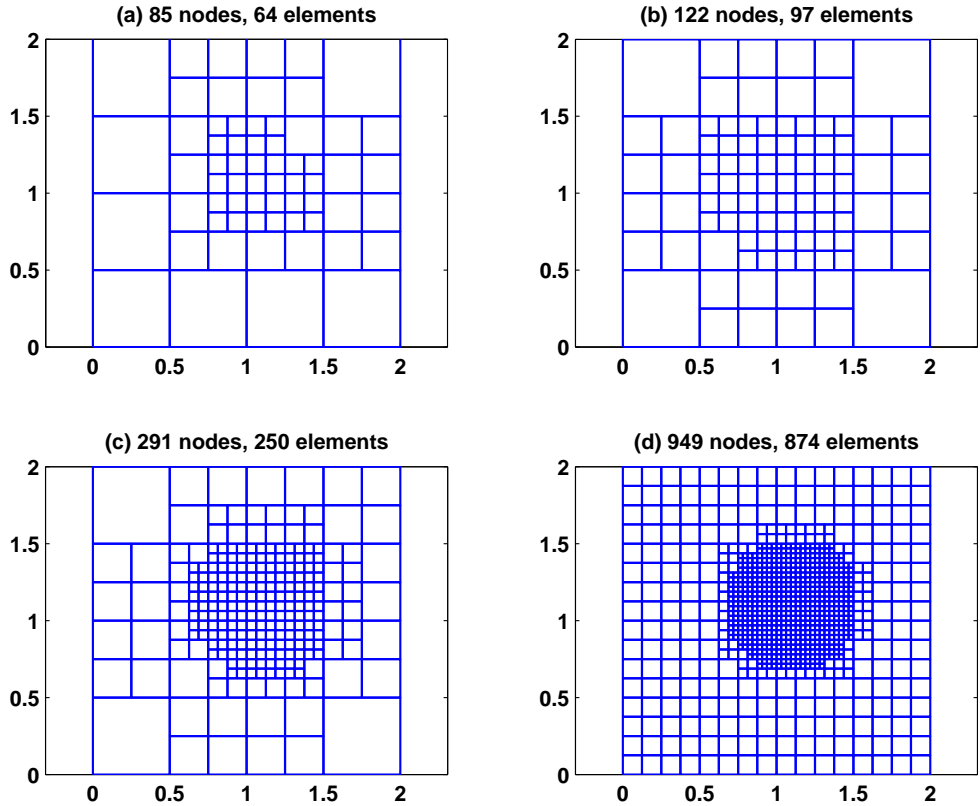}
\caption{Regularised quadtree meshes: 2-to-1 ratio is maintained.}
\label{fig:twotoonemesh}
\end{figure}

\begin{figure}[htpb]
\centering
\includegraphics[scale=0.8]{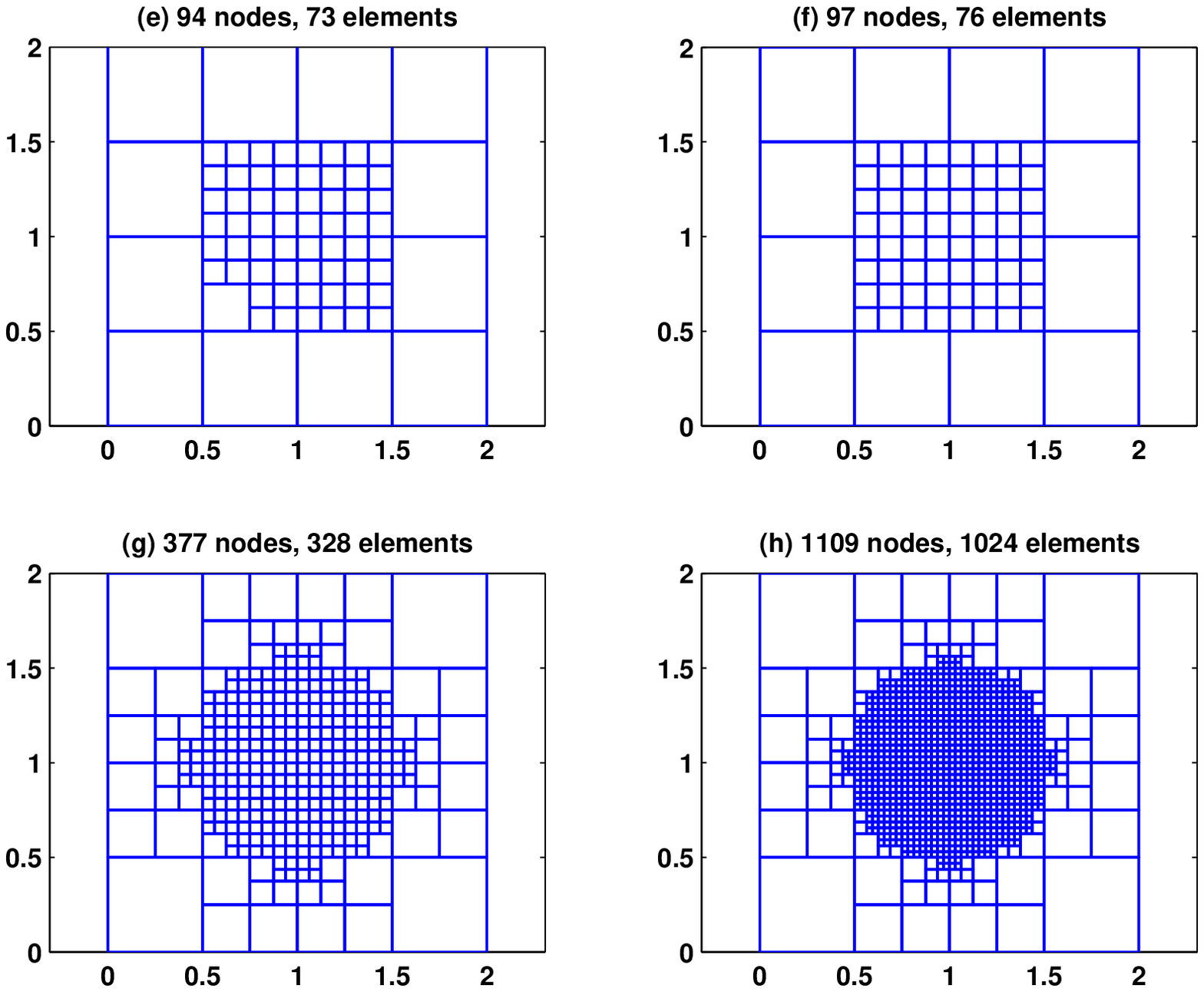}
\caption{Non-regularised quadtree meshes: 2-to-1 ratio is `not' maintained.}
\label{fig:notwotoonemesh}
\end{figure}

Consider the Laplace equation: $\nabla^2 u(\xx) = 0$ in $\Omega=$(0,1)$^2$. In this case, we consider two essential boundary condition: (a) Case A: Linear function, $g(\xx) = x+y$ and (b) Case B: Quadratic function, $g(\xx) = 1-x+5y-2xy-4x^2+4y^2$, imposed on $\partial \Omega$. The exact solution is: $u(\xx)=g(\xx)$. \frefs{fig:twotoonemesh} - (\ref{fig:notwotoonemesh}) show the quadtree mesh considered for this study. Two-to-one rule is maintained in \fref{fig:twotoonemesh} in which the elements have utmost only one hanging node along an edge, whilst in meshes shown in \fref{fig:notwotoonemesh}, some elements have more than one hanging node along an edge. The patch test results for various meshes and different approaches in the case of linear function is shown in Table \ref{table:laplinpatch}. It is seen that the FEM, the nSFEM and the SBFEM formulations over quadtree meshes pass patch test to machine precision. In the case of strain smoothing, both one subcell and $n$ subcells (triangular subcells) are considered. The error in the $L^2$ norm in the case of Laplace interplants is $\mathcal{O}($10$^{-8})$ when quadtree meshes with 2-to-1 rule is maintained and the accuracy slightly decreases when the element has more than one hanging node along its edges. This is consistent with the literature~\cite{tabarraeisukumar2005,tabarraeisukumar2007} and this can be attributed to the rational form of the Laplace shape functions. Moreover, the derivatives in the quadtree element are singular when using Laplace interpolant~\cite{tabarraeisukumar2005}. Also care must be taken to numerically integrate the terms when using Laplace interpolants.

\begin{table}[htpb]
\renewcommand\arraystretch{1.2}
\caption{Relative error in the $L^2$ norm for the patch test: Linear function.} 
\centering
\begin{tabular}{lrrrrrr}
\hline
Ratio & Mesh & Conforming shape& Laplace shape& \multicolumn{2}{c}{Strain Smoothing} & SBFEM \\
\cline{5-6}
 &  \frefs{fig:twotoonemesh} - (\ref{fig:notwotoonemesh}) & functions~\cite{gupta1978} & functions & 1 subcell & n subcell &  \\
 \hline
\multirow{3}{*}{2-to-1} & a & 4.68e$^{-16}$& 7.75e$^{-8}$ & 2.59e$^{-16}$ & 2.47e$^{-16}$ & 6.25e$^{-16}$ \\
& b & 1.65e$^{-15}$ & 6.62e$^{-8}$ & 2.22e$^{-15}$ & 1.79e$^{-15}$ & 1.05e$^{-15}$\\
& c & 5.50e$^{-15}$ & 4.26e$^{-8}$ & 5.67e$^{-15}$ & 8.16e$^{-15}$ & 8.02e$^{-15}$ \\
\cline{2-7}
\multirow{3}{*}{no 2-to-1} & e & - & 1.14e$^{-8}$ & 5.65e$^{-16}$ & 2.47e$^{-16}$ & 7.43e$^{-16}$ \\
& f & - & 2.47e$^{-7}$ & 8.12e$^{-16}$ & 5.32e$^{-15}$ & 6.66e$^{-16}$ \\
& g & - & 4.04e$^{-5}$ & 4.62e$^{-15}$ & 8.81e$^{-15}$ & 3.82e$^{-15}$ \\
& h & - & 7.64e$^{-7}$ & 1.54e$^{-14}$ & 1.69e$^{-14}$ & 1.70e$^{-14}$ \\
\hline
\end{tabular}
\label{table:laplinpatch}
\end{table}

\begin{table}[htpb]
\renewcommand\arraystretch{1.2}
\caption{Relative error in the $L^2$ norm for the patch test: quadratic function.} 
\centering
\begin{tabular}{lrrrrr}
\hline
Mesh & Conforming shape& Laplace shape& \multicolumn{2}{c}{Strain Smoothing} & SBFEM \\
\cline{4-5}
 \fref{fig:twotoonemesh} & functions~\cite{gupta1978} & functions & 1 subcell & n subcell & 1 \\
 \hline
 a & 3.39e$^{-3}$ & 4.38e$^{-3}$ & 1.46e$^{-2}$ & 3.25e$^{-3}$ & 4.75e$^{-3}$\\
 b & 2.19e$^{-3}$ & 2.63e$^{-3}$ & 1.66e$^{-2}$ & 2.25e$^{-3}$ & 2.93e$^{-3}$\\
 c & 9.18e$^{-4}$ & 1.15e$^{-3}$ & 8.10e$^{-3}$ & 1.30e$^{-3}$ & 1.37e$^{-3}$\\
 d & 1.81e$^{-4}$ & 2.11e$^{-4}$ & 9.80e$^{-4}$ & 1.88e$^{-4}$ & 2.25e$^{-4}$\\
 \hline
\end{tabular}
\label{table:lapquadpatch}
\end{table}

\subsection{Poisson problem} 
\begin{figure}[htpb]
\centering
\includegraphics[scale=0.8]{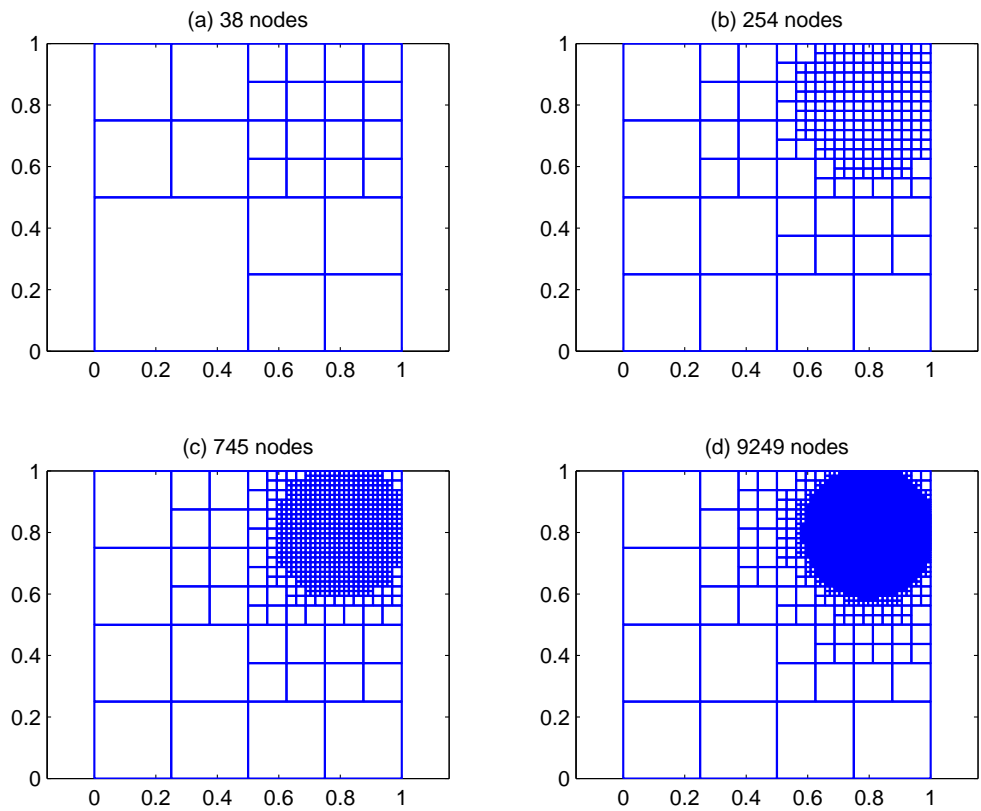}
\caption{Poisson problem: regularised quadtree meshes.}
\label{fig:poissonmesh}
\end{figure}

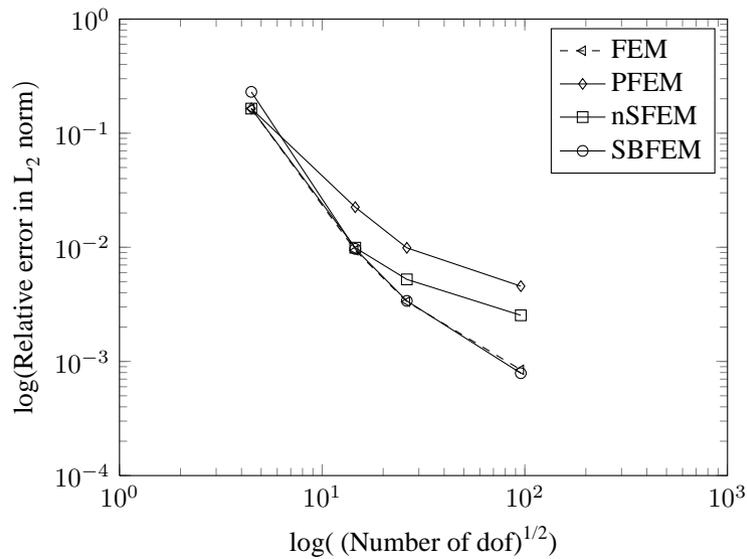
\begin{figure}
\centering
\newlength\figureheight 
\newlength\figurewidth 
\setlength\figureheight{6cm} 
\setlength\figurewidth{8cm}
%
%
%
%
\begin{tikzpicture}

\begin{axis}[%
width=\figurewidth,
height=\figureheight,
scale only axis,
xmode=log,
xmin=1,
xmax=1000,
xminorticks=true,
xlabel={$\text{log( (Number of dof)}^{\text{1/2}}\text{)}$},
ymode=log,
ymin=0.0001,
ymax=1,
yminorticks=true,
ylabel={$\text{log(Relative error in L}_\text{2}\text{ norm})$},
legend style={draw=black,fill=white,legend cell align=left}
]
\addplot [
color=black,
dashed,
mark=triangle,
mark options={solid,,rotate=90}
]
table[row sep=crcr]{
4.47213595499958 0.164135896618808\\
14.560219778561 0.00942483707832\\
26.1725046566048 0.003341938790257\\
95.5458005356593 0.0008447138759763\\
};
\addlegendentry{FEM};

\addplot [
color=black,
solid,
mark=diamond,
mark options={solid}
]
table[row sep=crcr]{
4.47213595499958 0.164330045731217\\
14.560219778561 0.02247540600683\\
26.1725046566048 0.00989330325124\\
95.5458005356593 0.004563412857985\\
};
\addlegendentry{PFEM};

\addplot [
color=black,
solid,
mark=square,
mark options={solid}
]
table[row sep=crcr]{
4.47213595499958 0.164103903579511\\
14.560219778561 0.009919676425528\\
26.1725046566048 0.005234305044628\\
95.5458005356593 0.002536500299614\\
};
\addlegendentry{nSFEM};

\addplot [
color=black,
solid,
mark=o,
mark options={solid}
]
table[row sep=crcr]{
4.47213595499958 0.230029595153746\\
14.560219778561 0.00968587614785\\
26.1725046566048 0.003390667313793\\
95.5458005356593 0.0007888666396518\\
};
\addlegendentry{SBFEM};

\end{axis}
\end{tikzpicture}%
\caption{Poisson problem: convergence in the relative error in the $L_2$ norm, where FEM - conventional FEM with conforming shape functions~\cite{gupta1978}, PFEM - polygonal FEM with Laplace interpolant~\cite{tabarraeisukumar2005}, nSFEM - cell-based smoothed FEM with $n$ triangular subcells and SBFEM - scaled boundary FEM.}
\label{fig:poissconve}
\end{figure}

In this example, we consider Poisson's equation over a square domain with Dirichlet boundary conditions imposed on the boundary. The governing equation and the boundary condition are:
\begin{align}
\nabla^2 u &= f, \hspace{0.2cm} \textup{in} \hspace{0.2cm} \Omega = (0,1)^2 \nonumber \\
u &= 0, \hspace{0.2cm} \textup{on} \hspace{0.2cm} \partial \Omega
\end{align}
The source term $f$ is chosen such that the exact solution of the problem is~\cite{tabarraeisukumar2005}: $u(\xx) = x^{10}y^{10}(1-x)(1-y)$. The quadtree mesh considered for this study are shown in \fref{fig:poissonmesh}. In this example, for the CSFEM, triangular subcells are considered. The mesh is refined where the steep gradient is anticipated. \fref{fig:poissconve} shows the convergence of the relative error in the $L_2$ norm for different approaches. It is seen that the FEM and the SBFEM yield similar results. The error from the PFEM formulation is higher when compared to other approaches. This could be attributed to the fact that the derivatives are singular at the hanging nodes and also to the accuracy of the numerical integration. All the approaches converge asymptotically with mesh refinement.

\subsection{Application to partition of unity method - weak discontinuity problem}
\begin{figure}
\centering
\subfigure[Geometry] {\scalebox{1}{\input{./Figures/circular.pstex_t}}}
\subfigure[Typical quadtree mesh]{\includegraphics[scale=0.5]{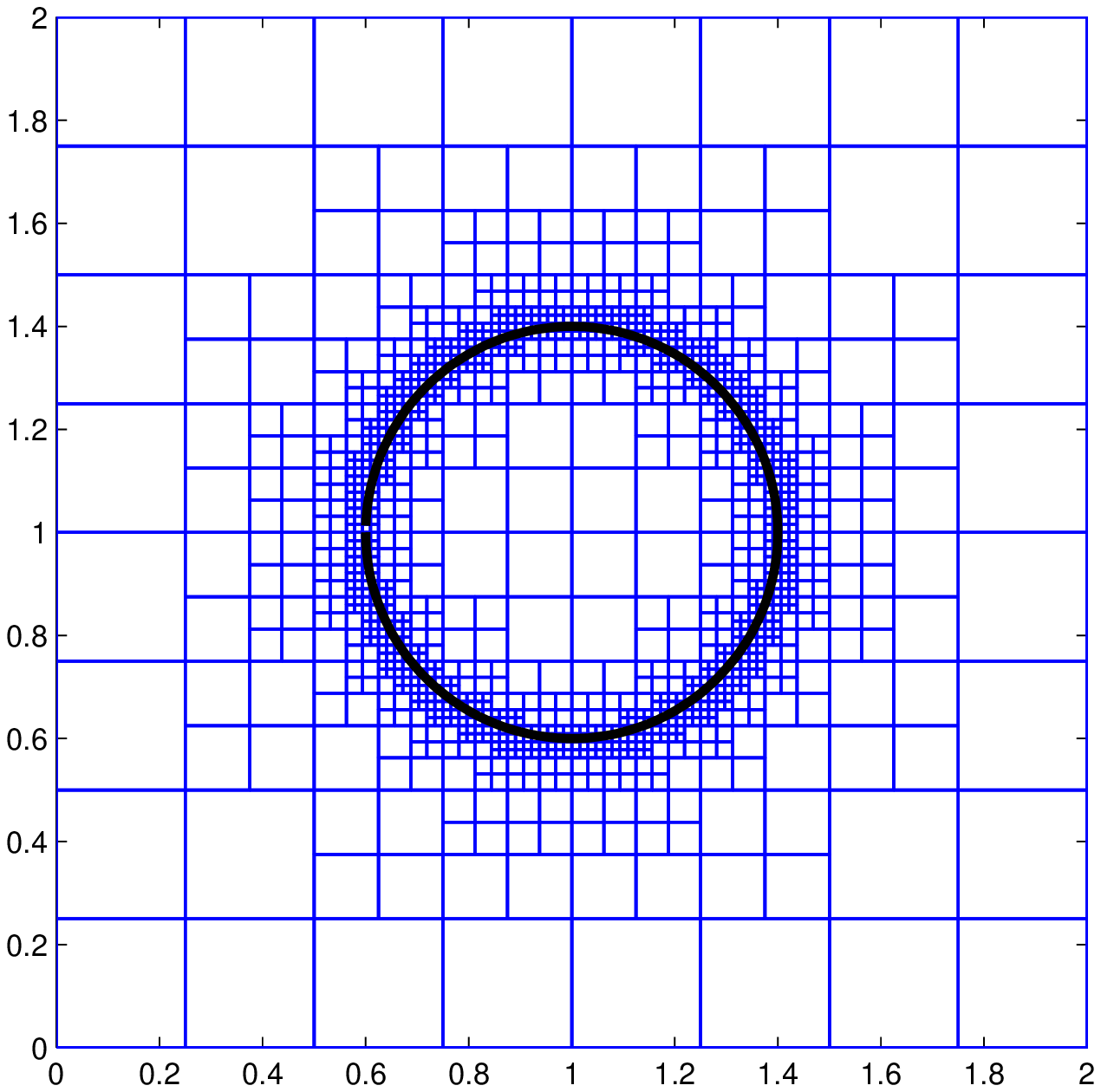}}
\caption{Bimaterial circular inhomogeneity boundary value problem: geometry and quadtree mesh (where `solid' black line represents the material interface). Note that the quadtree mesh does not conform to the material interface. Appropriate nodes are enriched whose nodal support is intersected by the material interface~\cite{sukumarchopp2001}.}
\label{fig:circular}
\end{figure}

In this example, the elastio-static response of a circular material inhomogeneity under radially symmetric loading as shown in \fref{fig:circular} is examined within the framework of the XFEM. One of the salient feature of the XFEM is that, by augmenting the FE approximation basis with additional functions, the local information of the problem can be retrieved without a need for a conforming mesh. As the information is local, a quadtree mesh is better suited than a uniformly refined mesh. \fref{fig:circular} shows a typical mesh used in this study. It is seen that the mesh is locally refined in the close proximity of the interface. A 2-to-1 ratio is maintained in this study. The material constants are constant within each domain, $\Omega_1$ and $\Omega_2$, but there is a material discontinuity across the interface, $\Gamma_1 (r=a)$. The Lam\'e constants in $\Omega_1$ and $\Omega_2$ are: $\lambda_1 = \mu_1 = 0.4$ and $\lambda_2 = 5.7692,~ \mu_2 = 3.8461$. These correspond to $E_1 = 1, ~\nu_1 = 0.25$ and $E_2 = 10, ~\nu_2 = 0.3$. A plane strain condition is assumed. A linear displacement field: $u_1 = x_1, ~u_2 = x_2 ~(u_r = r, ~u_{\theta} = 0)$ on the boundary $\Gamma_2 ~(r=b)$ is imposed and the governing equations and corresponding exact displacement solution is given in~\cite{sukumarchopp2001}.

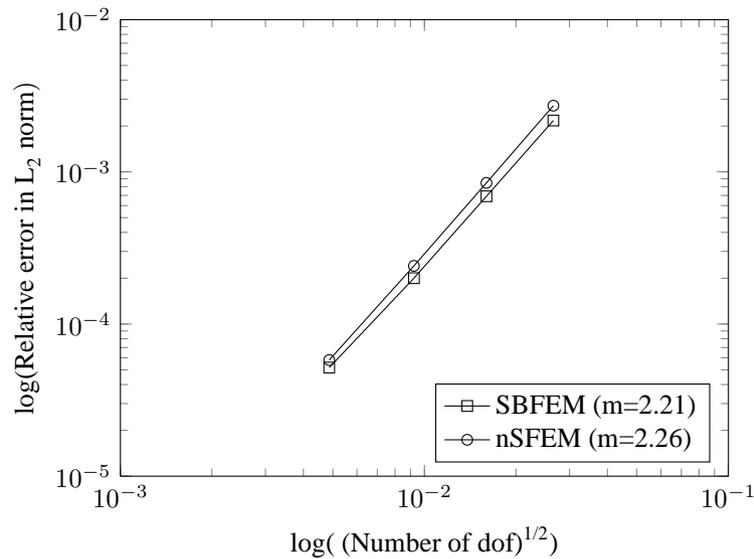
\begin{figure}
\centering
\setlength\figureheight{6cm} 
\setlength\figurewidth{8cm}
%
%
%
%
\begin{tikzpicture}

\begin{axis}[%
width=\figurewidth,
height=\figureheight,
scale only axis,
xmode=log,
xmin=0.001,
xmax=0.1,
xminorticks=true,
xlabel={$\text{log( (Number of dof)}^{\text{1/2}}\text{)}$},
ymode=log,
ymin=1e-05,
ymax=0.01,
yminorticks=true,
ylabel={$\text{log(Relative error in L}_\text{2}\text{ norm})$},
legend style={draw=black,fill=white,legend cell align=left},
legend pos = south east
]
\addplot [
color=black,
solid,
mark=square,
mark options={solid}
]
table[row sep=crcr]{
0.0265559525864072 0.002174094\\
0.0160005120245773 0.000691058482633321\\
0.00923632362561235 0.000200105412582593\\
0.00486320121702676 5.17e-05\\
};
\addlegendentry{SBFEM (m=2.21)};

\addplot [
color=black,
solid,
mark=o,
mark options={solid}
]
table[row sep=crcr]{
0.0265559525864072 0.002722104963245\\
0.0160005120245773 0.000846435133736582\\
0.00923632362561235 0.000240968478012676\\
0.00486320121702676 5.80150040529949e-05\\
};
\addlegendentry{nSFEM (m=2.26)};

\end{axis}
\end{tikzpicture}%
\caption{Bi-material circular inhomogeneity: the rate of convergence. The error is measured in displacement $L^2$ norm. $m$ is the average slope.}
\label{fig:circulardispnorm}
\end{figure}

For the present numerical study, a square domain of size $L \times L$ with $L=$2 is considered, where the outer radius is chosen to be $b=$2 and inner radius $a = $0.4. The rate of convergence of the relative error in the displacement $(L^2)$ norm is shown in \fref{fig:circulardispnorm}. In this example, for the CSFEM, triangular subcells are considered. It is seen that both the proposed approaches yield optimal convergence rate and with mesh refinement, both the approaches converge to analytical solution. The main advantage of the strain smoothing and the scaled boundary formulation is that sub-triangulation of the elements having hanging nodes is eliminated and explicit form of the shape function is not required. In this study, the hanging nodes are not enriched for simplicity, as hanging nodes with enrichment leads to additional computational difficulties. However, with the strain smoothing technique, the hanging nodes could be enriched and an approach described in~\cite{friesbyfut2011,legrainallais2011} can be adopted. For the sake of simplicity, this is not employed in this study.

\section{Conclusions}
In this paper, we presented two alternate approaches to treat hanging nodes in a quadtree mesh based on the cell-based strain smoothing and the scaled boundary polygonal formulation. The convergence and the accuracy of both approaches were demonstrated with numerical examples. It is seen that the approaches pass patch test and yield accurate results when compared with the polygonal formulation with Laplace interpolants. This may be attributed to the fact that the derivatives in a quadtree element are singular when employing the polygonal formulation~\cite{tabarraeisukumar2005}. Also a very high number of Gauss points (for example 25 points) are used in the polygonal formulation of an element having hanging nodes~\cite{tabarraeisukumar2005}. The presented techniques do not suffer from such difficulties. Moreover, the techniques do not require an explicit form of the shape function and no special numerical integration is required. In the case of scaled boundary polygon formulation, higher order element can be easily constructed along each edge. The application of the cell-based smoothing to three-dimensional problems is available in the literature~\cite{nguyen-xuannguyen2012}, whilst it is under investigation in the case of scaled boundary formulation. Nevertheless, the presented approaches can be extended to octree meshes in 3D which will be topic for the future communication.

\ack Sundararajan Natarajan would like to acknowledge the financial support of the School of Civil and Environmental Engineering, The University of New South Wales for his research fellowship since September 2012.

\bibliographystyle{wileyj}
\bibliography{qtree}
\end{document}